\documentclass[journal,comsoc]{IEEEtran}
\usepackage[T1]{fontenc}

\usepackage{cite}

\ifCLASSINFOpdf
\else
\fi
\usepackage{amsmath}
\usepackage{algorithmic}

\usepackage{array}
\usepackage{url}

\usepackage{algorithm}
\usepackage{algorithmic}
\usepackage{color}
\usepackage{amsmath}
\usepackage{url}
 \usepackage{amssymb}
\usepackage{amsfonts}
\usepackage{bm}
\usepackage{xcolor}
\usepackage{framed}
\usepackage{lipsum}
\usepackage{multirow}
\usepackage{threeparttable}
\usepackage{graphicx}
\usepackage{algorithmic}
\usepackage{mathrsfs}
\usepackage{array}
\usepackage{bm}
\usepackage{multirow}
\usepackage{booktabs}
\usepackage{nomencl}
\usepackage{dsfont}

\begin{document}
%
\title{Distributionally Robust Frequency Constrained Scheduling for an Integrated Electricity-Gas System}
%
%
%
\author{Lun~Yang,
        Yinliang~Xu,~\IEEEmembership{Senior Member,~IEEE,}
        Jianguo~Zhou,~\IEEEmembership{Member,~IEEE,}
        and~Hongbin~Sun,~\IEEEmembership{Fellow,~IEEE}
}

%
%

\markboth{Journal of \LaTeX\ Class Files,~Vol.~14, No.~8, August~2015}%
{Shell \MakeLowercase{\textit{et al.}}: Bare Demo of IEEEtran.cls for IEEE Communications Society Journals}
%



\maketitle

\begin{abstract}
Power systems are shifted from conventional bulk generation toward renewable generation. This trend leads to the frequency security problem due to the decline of system inertia. On the other hand, natural gas-fired units are frequently scheduled to provide operational flexibility due to their fast adjustment ability. The interdependence between power and natural gas systems is thus intensified. In this paper, we study the frequency constrained scheduling problem from the perspective of an integrated electricity-gas system under variable wind power. We propose a distributionally robust (DR) chance constrained optimization model to co-optimize the unit commitment and virtual inertia provision from wind farm systems. This model incorporates both frequency constraints and natural gas system (NGS) operational constraints and addresses the wind power uncertainty by designing DR joint chance constraints. We show that this model admits a mixed-integer second-order cone programming. Case studies demonstrate that the proposed approach can provide a highly reliable and computationally efficient solution and show the importance of incorporating NGS operational constraints in the frequency constrained scheduling problem.
\end{abstract}


%
\IEEEpeerreviewmaketitle
\section*{Nomenclature}
\subsection{Indices and Sets}
\addcontentsline{toc}{section}{Nomenclature}
\begin{IEEEdescription}[\IEEEusemathlabelsep\IEEEsetlabelwidth{$V_1,V_2,V_3$}]
    \item [$t \in \mathcal{T}$]  Set of time periods.
    \item [$g \in \mathcal{G}$] Set of gas-fired units (GFUs).
    \item [$n \in \mathcal{N}$] Set of non-GFUs.
    \item [$i \in \mathcal{G} \cup \mathcal{N} $] Set of GFUs and non-GFUs.
    \item[$w\in\mathcal{W}$] Set of wind farm systems (WFSs).
    \item[$l\in\mathcal{L}_e$] Set of transmission lines.
    \item[$d_e\in\mathcal{D}_e$] Set of electricity loads.
    \item[$m/n\in\Omega$] Set of gas nodes.
    \item[$s\in\mathcal{S}$] Set of gas sources.
    \item[$d_g \in \mathcal{D}_g$] Set of gas loads.
    \item[$ k\in \mathcal{C}$] Set of compressors.
    \item[$ (m,n)\in \mathcal{L}_g$] Set of pipelines.
    \item[$\mathcal{S}(m)$] Set of gas sources connected to node $m$.
    \item[$\mathcal{D}_g(m)$] Set of gas loads connected to node $m$.
    \item[$\mathcal{C}(m)$] Set of compressor inlet nodes connected to node $m$.
    \item[$\mathcal{G}(m)$] Set of GFUs connected to node $m$.

\end{IEEEdescription}
\subsection{Parameters}	
\addcontentsline{toc}{section}{Nomenclature}
\begin{IEEEdescription}[\IEEEusemathlabelsep\IEEEsetlabelwidth{$V_1,V_2,V_3$}]
    \item[$c_i, c_i^{B} $] Generation cost coefficients of generator $i$.
    \item[$c_i^{G},c_w^{W}$] PFR cost of generator $i$ and WFS $w$.
    \item[$c_i^{SU}, c_i^{SD}$] Start-up and shut-down cost of generator $i$.
    \item[$c_w^{V}$] Virtual inertia provision cost of WFS $w$.
    \item[$C_{mn}$] Weymouth constant of pipeline $m-n$.
    \item[$D_{d_e,t}$] Power load $d_e$ at hour $t$.
    \item[$F_{k}^C$] Maximum allowed gas flow of compressor $k$.
    \item[$F_l^{\text{max}}$] Capacity of transmission line $l$.
    \item[$\Delta f^{\text{DB}}$] Dead band of governors.
    \item[$H_i^G$] Inertia constant of generator $i$.
    \item[$H_w^W$] Virtual inertia constant of WFS $w$.
	\item[$P_t^D$] Load level at hour $t$.
    \item[$\Delta P_t^{Loss}$] Generation loss at hour $t$.
    \item[$\tilde{P}_{w,t}^W$] Wind power output of WFS $w$ at hour $t$.
    \item[$R_i^{G,\text{max}}$] Reserve limit of generator $i$.
    \item[$R_w^{W,\text{max}}$] Reserve limit of WFS $w$.
	\item[$T_d$] Delivery time.
    \item[$t_{\text{DB}}$] Dead-time band of governors.
    \item[$T_i^{on}, T_i^{off}$] Minimum-up/-down time of generator $i$.
    \item[$W_w^{\text{max}}$] Installed capacity of WFS $w$ .
    \item[$\Psi$] Shift-factor matrix.
    \item[$\varphi_g^G$] Conversion efficiency of GFU $g$.
     \item[$\nu_k^C$] Gas consumption percentage of compressor $k$.
\end{IEEEdescription}

\vspace{-0.4cm}
\subsection{Variables}	
\addcontentsline{toc}{section}{Nomenclature}
\begin{IEEEdescription}[\IEEEusemathlabelsep\IEEEsetlabelwidth{$V_1,V_2,V_3$}]
    \item[$F_{k,t}^C$] Gas flow through compressor $k$ at hour $t$.
    \item[$F_{d_g,t}^D$] Gas load $d_g$ at hour $t$.
    \item[$F_{g,t}^G$] Gas consumption of GFU $g$ at hour $t$.
    \item[$F_{s,t}^S$] Output of gas source $s$ at hour $t$.
    \item[$F_{mn,t}$] Gas flow within pipeline $m-n$ at hour $t$.
    \item[$F_{mn,t}^{in}$] In-gas flow within pipeline $m-n$ at hour $t$.
    \item[$F_{mn,t}^{out}$] Out-gas flow within pipeline $m-n$ at hour $t$.
    \item[$H_t$] Total system inertia at hour $t$.
    \item[$L_{mn,t}$] Linepack of pipeline $m-n$ at hour $t$.
    \item[$P_{w,t}^W$] Dispatched wind power of WFS $w$ at hour $t$.
    \item[$R_{i,t}^G$] PFR from generator $i$ at hour $t$.
    \item[$R_{w,t}^W$] PFR from WFS $w$ at hour $t$.
    \item[$x_{i,t}$] Operation status of generator $i$ at hour $t$.
    \item[$y_{w,t}$] Operation status of virtual inertia provision from WFS $w$ at hour $t$.
    \item[$z_{i,t}^u, z_{i,t}^d$] Start-up status and shut-down status of generator $i$ at hour $t$.
    \item[$\tau_{k,t}^C$] Gas consumption of compressor $k$ at hour $t$.
    \item[$\pi_{m,t}$] Pressure of node $m$ at hour $t$ .
    \item[$\pi_{k,t}^{in}, \pi_{k,t}^{out}$] Inlet-/outlet-node pressures of compressor $k$ at hour $t$.

\end{IEEEdescription}

\section{Introduction}
%
%
%
%
\IEEEPARstart{T}{he} increasing share of renewable energy such as wind generation is changing the structure of generation mix. With conventional thermal units being replaced by renewable generation units, the total system inertia will gradually decrease since wind generation units  are usually asynchronously interfaced to the grid via power electronic devices \cite{Doherty-TPS-2010}. Such trend raises the concern on frequency stability problems in case of a sudden imbalance between generation and demand \cite{Xie-PIEEE-2011}. In fact, some power systems due to the low inertia already face \emph{Rate of Change of Frequency} (RoCoF) issues, such as power systems in Texas \cite{Matevosyan-Springer-2017}, Ireland \cite{O'Sullivan-TPS-2014}, and northwest of China \cite{Wen-TPS-2016}.

To improve the frequency stability in the power system with high penetration of wind generation, researchers have revealed the importance of including frequency constraints (i.e., limitations on RoCoF, frequency nadir, and quasi-steady-state frequency) in power system scheduling models.
Restrepo \emph{et al.} \cite{Restrepo-TPS-2005} contribute to a mixed-integer linear programming (MILP) model for a deterministic unit commitment (UC) problem with primary frequency regulation constraints, in which only quasi-steady-state frequency limits are considered. Chang \emph{et al.} \cite{Chang-TPS-2013} estimate the frequency nadir based on historical data and include it into a UC problem. Ahmadi \emph{et al.}\cite{Ahmadi-TPS-2014} derive the analytical nonlinear frequency nadir constraints and incorporate them into a deterministic security-constrained UC problem by using the piecewise linearization technique. Zhang \emph{et al.} \cite{Zhang-TPS-2020} propose the concept of frequency security margin and incorporate it into the UC model, where the nonlinear frequency nadir constraint is also linearized by piecewise linear techniques. Trovato \emph{et al.} \cite{Trovato-TPS-2019} propose an MILP-based UC formulation that co-optimizes energy production and the multi-speed allocation of multi-speed frequency response services.

However, the above frequency constrained UC models are deterministic. To address the uncertainty in the low-inertia power system, Lee \emph{et al.} \cite{Lee-TPS-2013} present a two-stage stochastic frequency constrained economic dispatch model considering uncertainty from unit outage and adopt L-shape method to solve it. Teng \emph{et al.} \cite{Teng-TPS-2016} develop
a stochastic UC model that derives an MILP formulation for the frequency constraints. In this model, wind power uncertainty is described by a set of representative scenarios. Paturet \emph{et al.} \cite{Paturet-TPS-2020} propose an efficient piecewise piecewise linearization method to handle the nonlinear frequency nadir constraints in stochastic UC model and also adopt the scenario-based method to describe the wind power uncertainty and unit outage. The authors of \cite{Wen-TPS-2016}, \cite{Prakash-TPS-2018} employ interval optimization to model wind power uncertainty within the frequency constrained UC problem. More recently, Ding \emph {et al.} \cite{Ding-TPS-2021} consider the impact of uncertain wind generation on the virtual inertia provision from wind-storage systems and propose a two-stage chance constrained model to co-optimize the UC and virtual inertia. Then, the sample average approximation (SAA) is used to approximate the chance constraints of frequency limitations as MILP constraints.

The aforementioned studies are dedicated to discussing how to include frequency constraints into the UC model for the low-inertia power system and have made significant contributions to enhance the frequency dynamic performance. However, we notice that there are some important research gaps in existing studies, which are elaborated as below.

First, the aforementioned studies are conducted from the perspective of power systems. In many countries, natural gas is replacing coal as the predominant fuel for traditional bulk generation since gas-fired units (GFUs) hold the merits of lower carbon emission and higher operational flexibility relative to coal-fired units \cite{Wen-TSG-2018}. The deployment of GFUs enhances the interdependence between power and natural gas systems, which raises the necessity on coordinated analysis for both systems. Some inspiring works have been conduced in the following aspects: coordinated operation \cite{Chen-TPS-2020OE}-\cite{Wang-TSE-2018}, planning \cite{Zhao-TPS-2018}, and market \cite{Chen-TPS-2020}. However, to the best of the authors' knowledge, there is few work that investigates the impact of natural gas system (NGS) operational constraints on the frequency constrained scheduling of the low-inertia power system. In reality, state variations of NGS such as the shortage of gas sources and pipeline outage will directly impact the operation states of GFUs (e.g., start up and shut down) and then system inertia level. Motivated by this, this paper seeks to fill this gap and will study the frequency constrained scheduling problem for an integrated electricity-gas system (IEGS).

Second,
the existing studies above mainly employ scenario-based stochastic programming \cite{Teng-TPS-2016}, \cite{Paturet-TPS-2020}, interval optimization \cite{Wen-TPS-2016}, \cite{Prakash-TPS-2018}, and chance constrained programming \cite{Ding-TPS-2021} to deal with the wind power uncertainty. Scenario-based stochastic programming relies on the assumption on the particular probability distribution. However, the true distribution of uncertainty could be unknown in many practical cases. Moreover, scenario-based stochastic programming usually needs to generate a relatively large number of scenarios and thus leads to a high computational burden. Interval optimization describes the uncertainty within a certain range around a central forecast and just requires three non-probabilistic scenarios. This method will be more computational efficient compared to the scenario-based stochastic programming but may result in an overly conservative solution. As for chance constrained programming in \cite{Ding-TPS-2021}, the adopted SAA for approximating chance constraints as MILP constraints requires numerous samples to guarantee the performance of solutions and leads to a heavy computational effort due to the increase of introduced binary variables.

Distributionally robust (DR) chance constrained approach does not require the presumed probability distribution and large number of samples and characterizes the uncertainty by a family of underlying probability distributions, termed as ambiguity set \cite{Zheng-TSG-2021}. Recently, DR chance constrained approaches have been successfully applied to handle the uncertainty in optimal power flow \cite{Zhang-TPS-2017}, distribution system planning \cite{Zare-TPS-2018}, microgrid energy management \cite{Shi-TSG-2019}, and optimal power-gas flow (OPGF) problems \cite{Fang-APEN-2019}, \cite{Yang-TSG-2021}. It is worth mentioning that Chu \emph{et al.} \cite{Chu-TSG-2021} introduce a DR chance constrained approach to address the uncertain non-critical load shedding in the microgrid scheduling problem and formulate a DR individual chance constraint of frequency requirements under the moment-based ambiguity set. However, to the best of the authors' knowledge, the exploration of DR chance constrained approach in frequency constrained UC problem under wind power uncertainty has not been fully investigated.

Given the research gaps discussed above, we propose a DR frequency-constrained UC (DR-FCUC) method for an IEGS. This work is related somewhat to our previous work \cite{Yang-TSG-2021}, which proposes a DR chance constrained OPGF method. However, this work focuses on the frequency constrained scheduling problem and includes frequency constraints and UC constraints, which do not appear in \cite{Yang-TSG-2021}. Furthermore, to address wind power uncertainty, the work of \cite{Yang-TSG-2021} considers the DR individual chance constraints under the moment-based ambiguity set while this work designs the DR joint chance constraint under the moment-based ambiguity set, which is further extended to include unimodality information.

The main contributions of this work are summarized as follows:

1) We for the first time propose a DR-FCUC model considering IEGS operational constraints from the perspective of an IEGS to co-optimize the scheduling and virtual inertia from wind farm systems. The proposed DR-FCUC model addresses the wind power uncertainty by designing DR joint chance constraints under the moment-based ambiguity set. This is different from the DR individual chance constraint for handling the uncertain non-critical load shedding in \cite{Chu-TSG-2021}. DR joint chance constraint modeling gives a stronger guarantee on the solution reliability relative to the DR individual chance constraint modeling.

2) To obtain a tractable formulation for the proposed DR-FCUC model, we derive second-order cone (SOC) constraints for the DR joint chance constraints under the moment-based ambiguity set. We then adopt the penalty convex-concave procedure in \cite{Wang-TSE-2018}, \cite{Yang-TSG-2021} to reformulate the non-convex Weymouth equations as SOC constraints. The resulting model is a mixed-integer SOC programming (MISOCP), which can be easily solved by off-of-shelf solvers.

3) We extend the proposed DR-FCUC model to include the unimodality information to strengthen the moment-based ambiguity set and derive SOC constraints for the DR joint chance constraint under the ambiguity set with both moment and unimodality information. Incorporating unimodality information can reduce the conservativeness of the proposed DR-FCUC model under the moment-based ambiguity set.

\section{Modeling of frequency requirements with virtual inertia provision from WFSs }
In this section, we introduce the modeling of frequency requirements, namely limitations on RoCoF, frequency nadir, and quasi-steady-state frequency. As explained in \cite{Ding-TPS-2021}, wind farm system (WFS) controlled by power electronics can provide virtual inertia for the power system. Specifically, WFS can immediately provide extra power in response to the disturbance. Therefore, we consider the virtual inertia provision from WFSs in the modeling of frequency requirements.


As reported in \cite{Teng-TPS-2016}, constraints on the three frequency requirements can be derived from the the swing equation that describes the dynamics of system frequency deviation after a disturbance (e.g., generation loss, a sudden increased load).
The swing equation considering virtual inertia provision from WFSs is \cite{Ding-TPS-2021}:
\begin{align}
\label{SWING}
\nonumber 2H_t\frac{d\triangle f(\tau)}{d\tau}+D \cdot P_t^D\triangle f(\tau)=\sum_{i \in \mathcal {G}\cup \mathcal {N}} \triangle P_i(\tau)\\+\sum_{w \in \mathcal {W}} \triangle W_w(\tau)-\triangle P_t^{Loss}
\end{align}
where $\triangle f(\tau)$ is a frequency deviation after a contingency, $H_t$ is the system inertia, $D$ and $P_t^D$ are load damping rate and load level; $\triangle P_i(\tau)$ and $\triangle W_i(\tau)$ represent the power adjustment from synchronous generators (non-GFUs and GFUs) and WFSs following the generation loss $\triangle P_t^{Loss}$, respectively.

Following \cite{Ding-TPS-2021}, $\triangle P_i(\tau)$ and $\triangle W_i(\tau)$  in (\ref{SWING}) are represented by
\begin{align}
\label{PFRgen}
\triangle P_i(\tau)=\left\{\begin{array}{*{20}{l}}
0,~~~~~~~~~~~~~~~ \textrm{if}~~\tau\leq t_{\textrm{DB}}\\
\frac{R_{i,t}^G}{T_d}(\tau-t_{\text{DB}}),~\textrm{if}~~t_{\text{DB}}< \tau<T_d+t_{\text{DB}}\\
R_{i,t}^G,~~~~~~~~~~~~~\textrm{if}~~ \tau \geq T_d+t_{\text{DB}}
\end{array}\right.
\end{align}

\begin{align}
\label{PFRwind}
\triangle W_w(\tau)=\left\{\begin{array}{*{20}{l}}
0,~~~~~~~~~~~~~~~ \textrm{if}~~t\leq \tau_{\textrm{DB}}\\
\frac{R_{w,t}^W}{T_d}(\tau-t_{\text{DB}}),\textrm{if}~~t_{\text{DB}}< \tau<T_d+t_{\text{DB}}\\
R_{w,t}^W,~~~~~~~~~~~~\textrm{if}~~ \tau \geq T_d+t_{\text{DB}}
\end{array}\right.
\end{align}
where $t_{\text{DB}}$ and $T_d$ are dead-time band of governors and delivery time of frequency response, $R_{i,t}^G$ and $R_{w,t}^W$ are the primary frequency responses (PFRs) from synchronous units (non-GFUs and GFUs) and WFSs, respectively.

\subsubsection{RoCoF limit} In the short interval, it is reported that the post-contingency RoCoF is proportional to the amount of power shortage and inversely proportional to the system inertia \cite{Teng-TPS-2016}, \cite{Wen-TPS-2016}. The total system inertia $H_t$ that is required to satisfy the given requirement on the maximum RoCoF is modeled as
\begin{align}
\label{Rocof}
 \nonumber H_t=&\frac{\sum_{i\in \mathcal{G}\cup \mathcal{N}}H_i^G\cdot P_i^{\text{max}}\cdot x_{i,t}+\sum_{w\in \mathcal{W}}H_w^W\cdot W_w^{\text{max}}\cdot y_{w,t}}{f^0}\\
 &\geq \frac{|\triangle P_t^{Loss}|}{2RoCoF^{\text{max}}},~\forall t\in \mathcal{T}
\end{align}
where $f^0$ is the nominal frequency, $H_i^G$ and $H_w^W$ are the inertia constants of generator $i$ and WFS $w$, respectively. The binary variables $x_{i,t}$ and $y_{w,t}$ describe the operation status of generator $i$ and the operation status of virtual inertia provision from WFS $w$ at hour $t$.
\subsubsection{Frequency nadir limit} Given the requirement on the predefined threshold $\triangle f^{\text{max}}$,  $|\triangle f_{\text{nadir}}|$ is imposed to not exceed $\triangle f^{\text{max}}$. The frequency at nadir $|\triangle f_{\text{nadir}}|$ can be calculated by $\frac{|d\triangle f(t)|}{dt}=0$, giving rise to
\begin{align}
\label{fnadir}
\nonumber |\triangle f_{\text{nadir}}|&=\frac{2R_tH_t}{T_dD'^2_t}\text{log}\left(\frac{2R_t H_t}{T_dD'(\triangle P_i^{Loss}-D'_t\triangle f_{\text{DB}})+2R_t H_t}\right)\\
&+\frac{\triangle P_i^{Loss}-D'_t\triangle f_{\text{DB}}}{D'_t}+\triangle f_{\text{DB}} \leq \triangle f^{\text{max}}
\end{align}
where $R_t=\sum_{i\in \mathcal{G} \cup \mathcal{N}}R_{i,t}^G+\sum_{w\in \mathcal{W}}R_{w,t}^W$ represents the total PFR at hour $t$, and $D'_t=D\cdot P^D_t$.

Constraint (\ref{fnadir}) with a complex form is hard to be directly applied to the UC problem. The sufficient condition for constraint (\ref{fnadir}) to be satisfied is enforcing \cite{Teng-TPS-2016}:
\begin{align}
\label{Fnadirlimit}
 \nonumber R_tH_t&=R_t\frac{\sum_{i\in \mathcal{G}\cup\mathcal{N}}H_i^G\cdot P_i^{\text{max}}\cdot x_{i,t}+\sum_{w\in \mathcal{W}}H_w^W\cdot W_w^{\text{max}}\cdot y_{w,t}}{f^0}\\
 &\geq \kappa_t, \forall t\in \mathcal{T}
\end{align}
where $\kappa_t$ is the unique solution obtained from
\begin{align}
\label{fnadirsolution}
\nonumber &\frac{2\kappa_t}{T_d}\text{log}\left(\frac{2\kappa_t}{T_dD'(\triangle P_t^{Loss}-D'\triangle f_{\text{DB}}) +2\kappa_t }\right)\\
&=D'^2\left( \triangle f^{\text{max}}-\triangle f_{\text{DB}} \right)-D'\left(\triangle P_t^{Loss}-D'\triangle f_{\text{DB}}\right).
\end{align}
The proof of the sufficient condition can be found in \cite{Teng-TPS-2016}.

Note that constraint (\ref{Fnadirlimit}) is bilinear due to $R_t\cdot x_{i,t}$ and $R_t\cdot y_{w,t}$. Since the form of bilinear term is the product of a continuous variable and a binary variable, constraint (\ref{Fnadirlimit}) can be exactly reformulated as the following MILP constraints:
\begin{subequations}
\begin{align}
\label{fnadirMILP}
\nonumber \frac{\sum_{i\in \mathcal{G}\cup \mathcal{N}}H_i^g\cdot P_i^{\text{max}}\cdot X_{i,t}+\sum_{w\in \mathcal{W}}H_w^W\cdot W_w^{\text{max}}\cdot Y_{w,t}}{f^0}
\\\geq \kappa_t,\forall t \in \mathcal{T}
\end{align}
\begin{align}
\label{BigMgenOFF}
-Mx_{i,t}\leq X_{i,t} \leq Mx_{i,t},\forall i\in \mathcal{G}\cup  \mathcal{N}, t \in \mathcal{T}
\end{align}
\begin{align}
\label{BigMgenON}
\nonumber -M(1-x_{i,t})&\leq X_{i,t}-\sum_{i\in \mathcal{N}\cup\mathcal{G}}R_{i,t}^G -\sum_{w\in \mathcal{W}}R_{w,t}^W
\\&\leq M(1-x_{i,t}),\forall i\in \mathcal{G}\cup \mathcal{N} , t \in \mathcal{T}
\end{align}
\begin{align}
\label{BigMwindOFF}
-My_{w,t}\leq Y_{w,t} \leq My_{w,t},\forall w\in \mathcal{W}, t \in \mathcal{T}
\end{align}
\begin{align}
\label{BigMwindON}
\nonumber -M(1-y_{w,t})\leq Y_{w,t}-\sum_{i\in \mathcal{N}\cup \mathcal{G}}R_{i,t}^G -\sum_{w\in \mathcal{W}}R_{w,t}^W \\\leq M(1-y_{w,t}),\forall w\in \mathcal{W}, t \in \mathcal{T}
\end{align}
\end{subequations}
where $M$ is a predefined big-positive number, $X_{i,t}$ and $Y_{w,t}$ are introduced auxiliary variables.

\subsubsection{Quasi-steady-state frequency limit} Given the allowable maximum quasi steady-state frequency deviation $\triangle f_{\text{qss}}^{\text{max}}$, the quasi-steady-state frequency limit is modeled by
\begin{align}
\label{QSSlimit}
\frac{\triangle P_i^{Loss}-R_t}{D'} \leq \triangle f_{\text{qss}}^{\text{max}} , t \in \mathcal{T}.
\end{align}

\section{Frequency constrained DR-UC formulation for an IEGS}
In this section, we first describe the UC formulation considering NGS operational constraints for an IEGS and then introduce the DR joint chance constraint modeling for wind power uncertainty in Section III.E. Finally, we include the frequency constraints in the UC formulation and formulate the overall DR frequency-constrained UC (DR-FCUC) problem for an IEGS in Section III.F.
\vspace{-0.3cm}
\subsection{Objective function}
Following \cite{Ding-TPS-2021}, the objective function of the DR-FCUC problem is to minimize the total cost:
\begin{align}
\label{obj}
\nonumber & \min \sum_{t\in\mathcal{T}}\sum_{i\in \mathcal{G} \cup \mathcal{N}}(c_i^{SU}z_{i,t}^u+c_i^{SD}z_{i,t}^d+c_i^Bx_{i,t}
\\&+c_iP_{i,t}+c_i^GR_{i,t}^G)+\sum_{t\in\mathcal{T}}\sum_{w\in\mathcal{W}}\left(c_w^Vy_{w,t}+c_w^WR_{w,t}^W\right)
\end{align}
where the first five terms represent the start-up cost, shut-down cost, no-load cost, running cost, and PFR cost of GFUs and non-GFUs. The last two terms represent the cost of virtual inertia provision and PFR cost of WFSs.

\vspace{-0.5cm}
\subsection{Natural Gas System Operational Constraints}
We adopt a nonconvex gas flow model with linepack to describe the natural gas system \cite{Wang-TSE-2018}. The natural gas system operational constraints are introduced below.
\subsubsection{Nodal gas flow balance constraints}
\begin{align}
\label{GASbalance}
\nonumber &\sum_{s\in\mathcal{S}(m)} F_{s,t}^S-\sum_{g\in\mathcal{G}(m)}F_{g,t}^G-\sum_{d_g\in\mathcal{D}_g(m)}F_{d_g,t}^D-\sum_{k\in \mathcal{C}(m)}\tau_{k, t}^C\\
&=\sum_{n\in \Omega(m)}F_{mn,t}+\sum_{k\in \mathcal{C}(m)}F_{k,t}^C, \forall m \in \Omega, t\in \mathcal{T}
\end{align}
where (\ref{GASbalance}) imposes the nodal gas flow balance.
\subsubsection{Gas pipeline constraints}
\begin{subequations}
\begin{align}
\label{avegasflow}
&F_{mn,t}=\frac{F_{mn,t}^{in}+F_{mn,t}^{out}}{2}, \forall (m,n)\in \mathcal{L}_g, t\in\mathcal{T}\\
\label{Weymouth}
&F_{mn,t}^2=C_{mn}^2(\pi_{m,t}^2-\pi_{n,t}^2),\forall (m,n)\in \mathcal{L}_g, t\in\mathcal{T}\\
\label{Linepackpipe}
&LP_{mn,t}=K_{mn}\frac{\pi_{m,t}+\pi_{n,t}}{2}, \forall (m,n)\in \mathcal{L}_g, t\in\mathcal{T}\\
\label{Linepackrelation}
\nonumber &F_{mn,t}^{in}-F_{mn,t}^{out}=LP_{mn,t}-LP_{mn,t-1}\\
&\forall (m,n)\in \mathcal{L}_g, t\in\mathcal{T}\\
\label{Finallinepack}
&\sum_{(m,n)\in\mathcal{L}_g}LP_{mn,T} \geq \sum_{(m,n)\in\mathcal{L}_g}LP_{mn,0}
\end{align}
\end{subequations}
where (\ref{avegasflow}) computes the average gas flow within a pipeline. Constraint (\ref{Weymouth}) is the Weymouth equation that describes the relationship between average gas flow and two-end pressures of the pipeline. Note that there is an implicit assumption used in (\ref{Weymouth}), i.e., $F_{mn,t}\geq 0$, which means that the gas flow directions are known \emph{a priori}. This assumption is reasonable and commonly-used in day-ahead electricity-gas coordinated operations considering the fact that the gas flow directions do not alter intra-day \cite{Wang-TSE-2018}, \cite{Chen-TPS-2020}. Constraint (\ref{Linepackpipe}) defines the linepack within a pipeline based on the average pressures of two ends. Constraint (\ref{Linepackrelation})  establishes the relationship between the linepack and in/out gas flows in a pipeline. Constraint (\ref{Finallinepack})  describes that the final-period total linepack needs to be restored for its usage at the next horizon \cite{Fang-TSE-2018}.

\subsubsection{Operational constraints of compressors}
\begin{subequations}
\label{compressors}
\begin{align}
\label{Comflow}
&0\leq F_{k,t}^C\leq F_{k}^{C, \text{max}}, \forall k\in \mathcal{C}, t\in\mathcal{T}\\
\label{Comloss}
&\tau_{k,t}^C=\nu_{k}^CF_{k,t}^C, \forall k\in \mathcal{C}, t\in\mathcal{T}\\
\label{Comratio}
&R_{k}^{C,\text{min}}\pi_{k,t}^{in}\leq \pi_{k,t}^{out}\leq R_{k}^{C,\text{max}}\pi_{k,t}^{in}, \forall k\in \mathcal{C}, t\in\mathcal{T}
\end{align}
\end{subequations}
where (\ref{Comflow}) restricts the transported gas flow of each compressor. Constraint (\ref{Comloss}) computes gas consumption caused by gas-driven compressors \cite{Wang-TSE-2018}. Constraint (\ref{Comratio}) limits the compressor ratio.
\subsubsection{Gas source constraints}
\begin{align}
\label{GASsources}
F_{s,t}^{\text{min}}\leq F_{s,t}\leq F_{s,t}^{\text{max}}, \forall s\in \mathcal{S}, t\in\mathcal{T}
\end{align}
where (\ref{GASsources}) restricts the output of each gas source.
\subsubsection{Nodal pressure limits}
\begin{align}
\label{GASpress}
\pi_m^{\text{min}}\leq \pi_{m,t}\leq \pi_m^{\text{max}},\forall m\in\Omega, t\in\mathcal{T}
\end{align}
where (\ref{GASpress}) restricts the nodal pressure.

\vspace{-0.5cm}
\subsection{Power System Operational Constraints}
Similar to \cite{Teng-TPS-2016}-\cite{Prakash-TPS-2018}, we use a shift-based DC power flow model to describe the power system. The power system operational constraints are introduced below.
\subsubsection{Operational constraints of generators}
\begin{subequations}
\label{generators}
\begin{align}
\label{logicrelation}
&z_{i,t}^u-z_{i,t}^d=x_{i,t}-x_{i,t-1},\forall i \in \mathcal{G} \cup \mathcal{N}, t\in \mathcal{T}\\
\label{ONOFFBinary}
&z_{i,t}^u, z_{i,t}^d, x_{i,t} \in \{0,1\},\forall i \in \mathcal{G} \cup \mathcal{N}, t\in \mathcal{T}\\
\label{minON}
&\sum_{\tau=\max\{1,t-T_i^{on}+1\}}^t z_{i,t}^u \leq x_{i,t},\forall i \in \mathcal{G} \cup \mathcal{N}, t\in \mathcal{T}\\
\label{minOFF}
&\sum_{\tau=\max\{1,t-T_i^{off}+1\}}^t z_{i,t}^d \leq 1-x_{i,t},\forall i \in \mathcal{G} \cup \mathcal{N}, t\in \mathcal{T}\\
\label{geneLIMIT}
&x_{i,t}P_{i}^{\min}\leq P_{i,t}\leq x_{i,t}P_{i}^{\max}-R_{i,t}^G,\forall i \in \mathcal{G} \cup \mathcal{N}, t\in \mathcal{T}\\
\label{Rgene}
&0\leq R_{i,t}^G\leq x_{i,t}R_{i}^{G,\max},\forall i \in \mathcal{G} \cup \mathcal{N}, t\in \mathcal{T}\\
\label{Ramping}
\nonumber &-RD_i\leq P_{i,t}+ R_{i,t}^G-P_{i,t-1}- R_{i,t-1}^G\leq RU_i,\\&\forall i \in \mathcal{G} \cup \mathcal{N}, t\in \mathcal{T}
\end{align}
\end{subequations}
where (\ref{logicrelation}) and (\ref{ONOFFBinary}) describe the logic relationship of generators status. Constraints (\ref{minON}) and (\ref{minOFF}) state the minimum-up and minimum-down time limits. Constraint (\ref{geneLIMIT}) ensures that the
generation output in combination with reserve capacity is within the generation limits. Constraint (\ref{Rgene}) ensures that the PFR is within the reserve capacity limits. Constraint (\ref{Ramping}) describes the ramping limits of each generator.
\subsubsection{Operational constraints of WFSs}
\begin{subequations}
\begin{align}
\label{windLIMIT}
& P_{w,t}^W+R_{w,t}^W \leq \tilde{P}_{w,t}^{W},\forall w \in \mathcal{W}, t\in \mathcal{T}\\
\label{Rwind}
&0\leq R_{w,t}^W\leq y_{w,t}R_{w}^{W,\max},\forall w \in \mathcal{W}, t\in \mathcal{T}\\
\label{scheduledWIND}
&P_{w,t}^W \geq 0,\forall w \in \mathcal{W}, t\in \mathcal{T}
\end{align}
\end{subequations}
where (\ref{windLIMIT}) limits the scheduled wind power $P_{w,t}^W$ in combination with the PFR $R_{w,t}^W$ to be less than the available wind power output $\tilde{P}_{w,t}^{W}$. Constraint (\ref{Rwind}) describes the reserve capacity limits of each WFS. Constraint (\ref{scheduledWIND}) ensures that the scheduled wind power $P_{w,t}^W$ is nonnegative.

\subsubsection{Transmission capacity limits}
\begin{align}
\label{LINElimit}
&\nonumber -F_l^{\max}\leq \sum_{i \in \mathcal{G} \cup \mathcal{N}} \Psi_{l,i} P_{i,t}+\sum_{w\in\mathcal{W}}\Psi_{l,w} P_{w,t}^W\\&-\sum_{d_e \in \mathcal{D}_e} \Psi_{l,d_e} D_{d_e,t}
\leq F_l^{\max}, \forall l\in \mathcal{L}_e, t\in \mathcal{T}
\end{align}
where (\ref{LINElimit}) limits the power flow of line within its transmission capacity.
\subsubsection{Nodal power balance constraints}
\begin{align}
\label{Powerbalance}
\sum_{i \in \mathcal{G} \cup \mathcal{N}} P_{i,t}+\sum_{w \in \mathcal{W}} P_{w,t}^W=\sum_{d_e \in \mathcal{D}_e} D_{d_e,t},\forall t\in \mathcal{T}
\end{align}
where (\ref{Powerbalance}) ensures the system-wide power balance.

\vspace{-0.3cm}
\subsection{Coupling Constraint}
Power system and natural gas system are physically linked by GFUs, where the coupling constraint is modeled by \cite{Wang-TSE-2018}:
\begin{align}
\label{coupling}
F_{g,t}^G=\varphi_g^G (P_{g,t}+R_{g,t}^G),\forall g\in \mathcal{G}, t\in \mathcal{T}.
\end{align}
\vspace{-0.8cm}
\subsection{DR Joint Chance Constraints}
We note that constraint (\ref{windLIMIT}) involving wind power output $\tilde{P}_{j,t}^W$ does not consider the uncertainty associated with wind power.
In this paper, we consider a set of underlying probability distributions (termed as ambiguity set) to address the uncertainty and design the following DR joint chance constraint:
\begin{align}
\label{DRCC}
\inf_{\mathbb{P}\in \mathcal{D}}\mathbb{P}\left\{P_{w,t}^W+R_{w,t}^W \leq \tilde{P}_{w,t}^{W}, \forall w \in \mathcal{W}\right\}\geq 1-\epsilon,\forall t\in \mathcal{T}.
\end{align}
Constraint (\ref{DRCC}) enforces that for all wind farms at hour $t$, the scheduled wind power $P_{w,t}^W$ in combination with the PFR $R_{w,t}^W$ are simultaneously less than the wind power output $\tilde{P}_{w,t}^{W}$ with at least joint probability $1-\epsilon$ for all probability distributions of $\mathcal{D}$.  Similar to \cite{Zare-TPS-2018}-\cite{Fang-APEN-2019}, we consider a moment-based ambiguity set defined on mean and variance information to characterize the probability distributions of $\tilde{P}_{w,t}^W$:
\begin{align}
\label{ambiguityset}
\mathcal{D}=\left\{\mathbb{P}:
\begin{array}{*{20}{l}}
\mathbb{E}_\mathbb{P}[\tilde{P}_{w,t}^{W}]=\mu_{w,t}, \mathbb{E}_\mathbb{P}\left[(\tilde{P}_{w,t}^{W}-\mu_{w,t})^2\right]=\sigma_{w,t}^2, \\
\forall w \in \mathcal{W}, t\in \mathcal{T}, \\
\end{array}
\right\}
\end{align}
where $\mu_{w,t}$ and $\sigma_{w,t}^2 $ are mean and variance of $\tilde{P}_{w,t}^W$, respectively.

\subsection{DR Frequency Constrained UC Formulation}
Based on the contents in Section II, Section III.A-E, the overall DR frequency constrained UC (DR-FCUC) for an IEGS is formulated as
\begin{align}
\label{DRUC}
\begin{array}{*{20}{l}}
\min ~(\ref{obj})\\
\text{s.t.~(\ref{Rocof}), (\ref{fnadirMILP})-(\ref{BigMwindON}), (\ref{QSSlimit})}\\
\text{~~~~(\ref{GASbalance}), (\ref{avegasflow})-(\ref{Finallinepack}), (\ref{Comflow})-(\ref{Comratio}), (\ref{GASsources}), (\ref{GASpress}), (\ref{coupling})}\\
\text{~~~~(\ref{logicrelation})-(\ref{Ramping}), (\ref{Rwind}), (\ref{scheduledWIND}), (\ref{LINElimit}), (\ref{Powerbalance}), (\ref{DRCC})}
\end{array}
\end{align}
where the first line constraints (\ref{Rocof}), (\ref{fnadirMILP})-(\ref{BigMwindON}), (\ref{QSSlimit}) are frequency constraints (limitations on RoCoF, frequency nadir, and quasi-steady-state frequency). The second line constraints (\ref{GASbalance}), (\ref{avegasflow})-(\ref{Finallinepack}), (\ref{Comflow})-(\ref{Comratio}), (\ref{GASpress}), (\ref{coupling}) are NGS operational constraints. The third line constraints (\ref{logicrelation})-(\ref{Ramping}), (\ref{Rwind}), (\ref{scheduledWIND}), (\ref{LINElimit}), (\ref{Powerbalance}), (\ref{DRCC}) are power system operational constraints.
\vspace{-0.2cm}
\section{Solution Methodology}
The proposed DR-FCUC model (\ref{DRUC}) is computationally intractable due to the DR joint chance constraint (\ref{DRCC}), and the nonconvex Weymouth equation (\ref{Weymouth}). In this section, we introduce the corresponding treatments.
\vspace{-0.2cm}
\subsection{Convex reformulation for DR joint chance constraints}
We derive a convex inner approximation to reformulate the DR joint chance constraint (\ref{DRCC}) as SOC constraints.

\emph{Theorem 1: For all $t\in\mathcal{T}$, the DR joint chance constraint (\ref{DRCC}) under the ambiguity set (\ref{ambiguityset}) is inner approximated by the following SOC constraints:}
\begin{subequations}
\label{DRJCCtheorem}
\begin{align}
\label{CHBF}
&r_{w,t}\sigma_{w,t}\leq -P_{w,t}^W-R_{w,t}^W+\mu_{w,t},\forall w \in \mathcal{W}\\
\label{soc1}
&\left\|
\begin{array}{*{20}{l}}
2s_{w,t}\\
1
\end{array}\right\|_2 \leq 2\epsilon_{w,t}+1, \forall w \in\mathcal{W}\\
\label{soc2}
&\left\|
\begin{array}{*{20}{l}}
2\\
r_{w,t}-s_{w,t}
\end{array}\right\|_2 \leq r_{w,t}+s_{w,t}, \forall w\in\mathcal{W}\\
\label{BFsumT}
&\sum_{w\in\mathcal{W}}\epsilon_{w,t} \leq \epsilon\\
\label{Theoremvariable}
&\epsilon_{w,t}, r_{w,t}, s_{w,t} \geq 0, \forall w \in \mathcal{W}.
\end{align}
\end{subequations}

\emph{Proof:} See the Appendix A.
\vspace{-0.3cm}
\subsection{Convexification for non-convex Weymouth equations}
The nonconvex Weymouth equation (\ref{Weymouth}) is equivalently casted as two opposite constraints:
\begin{subequations}
\begin{align}
\label{WeymouthSOC}
&\frac{F_{mn,t}^2}{C_{mn}^2} \leq \pi_{m,t}^2-\pi_{n,t}^2, \forall (m,n)\in \mathcal{L}_g, t\in\mathcal{T}\\
\label{WeymouthCONCAVE}
&\pi_{m,t}^2\leq \frac{F_{mn,t}^2}{C_{mn}^2}+\pi_{n,t}^2, \forall (m,n)\in \mathcal{L}_g, t\in\mathcal{T}
\end{align}
\end{subequations}
where (\ref{WeymouthSOC}) is an SOC constraint, whose standard SOC form is
\begin{align}
\label{WeymouthSTANDSOC}
&\left\|
\begin{array}{*{20}{l}}
F_{mn,t}/ C_{mn}\\
\pi_{n,t}
\end{array}\right\|_2 \leq \pi_{m,t}, \forall (m,n)\in \mathcal{L}_g, t\in\mathcal{T}.
\end{align}

Constraint (\ref{WeymouthCONCAVE}) is concave and difference-of-convex program form. We introduce the penalty convex-concave procedure (PCCP) in \cite{Yang-TSG-2021}, \cite{Lipp-OE-2016} to address it.

Given a linearization point $\left(F_{mn,t}^r, \pi_{n,t}^r\right )$ obtained from last iteration $r$, the first-order Taylor series expansion of the right-hand side in (\ref{WeymouthCONCAVE}) is
\begin{align}
\label{Taylorapproxi}
\nonumber &\frac{F_{mn,t}^2}{C_{mn}^2}+\pi_{n,t}^2 \approx \frac{2F_{mn,t}^rF_{mn,t}-(F_{mn,t}^r)^2}{C_{mn}^2}\\&+2\pi_{n,t}^r\pi_{n,t}-(\pi_{n,t}^r)^2, \forall (m,n)\in \mathcal{L}_g, t\in\mathcal{T}.
\end{align}

Together with (\ref{Taylorapproxi}), constraint (\ref{WeymouthCONCAVE}) is approximated as \cite{Wang-TSE-2018}:
\begin{subequations}
\label{Taylor}
\begin{align}
\label{TaylorapproxiS}
\nonumber &\pi_{m,t}^2\leq  \frac{2F_{mn,t}^rF_{mn,t}-(F_{mn,t}^r)^2}{C_{mn}^2}\\&+2\pi_{n,t}^r\pi_{n,t}-(\pi_{n,t}^r)^2 +s_{mn,t}^r, \forall (m,n)\in \mathcal{L}_g, t\in\mathcal{T},\\
\label{slackVariable}
&s_{mn,t}^r \geq 0, \forall (m,n)\in \mathcal{L}_g, t\in\mathcal{T}
\end{align}
\end{subequations}
where $s_{mn,t}^r$ is a negative slack variable.

\begin{algorithm}[t]
 \caption{A sequential algorithm based on PCCP }
 \begin{algorithmic}[1]
\renewcommand{\algorithmicrequire}{\textbf{Initialization:}}
\REQUIRE Set the tolerance $\varepsilon_{gap}$, maximum iteration $R^{max}$, and parameters of $\varrho_0$, $\varrho^{max}$, and $\kappa$.
\STATE  Set the iteration index $r=0$. Solve the relaxed FCDRUC \textbf{P1} in (\ref{SOCDRUC}) to obtain the solutions $F_{mn,t}^r$ and $\pi_{n,t}^r$. Calculate the maximum SOC relaxation gap $ M_{gap}=\max\left\{(\pi_{m,t}^2-\pi_{n,t}^2-F_{mn,t}^2/C_{mn}^2)/\pi_{m,t}^2\right\}, \forall (m,n)\in\mathcal{L}_g, t\in\mathcal{T}.$
\begin{align}
\label{SOCDRUC}
\begin{array}{*{20}{l}}
\textbf{P1}: \min ~(\ref{obj})\\
\text{s.t.~(\ref{Rocof}), (\ref{fnadirMILP})-(\ref{BigMwindON}), (\ref{QSSlimit})}\\
\text{(\ref{GASbalance}), (\ref{avegasflow}), (\ref{Linepackpipe})-(\ref{Finallinepack}), (\ref{compressors}), (\ref{GASsources}) (\ref{GASpress}), (\ref{coupling}), (\ref{WeymouthSOC})}\\
\text{(\ref{logicrelation})-(\ref{Ramping}), (\ref{Rwind}), (\ref{scheduledWIND}), (\ref{LINElimit}), (\ref{Powerbalance}), (\ref{DRJCCtheorem})}.
\end{array}
\end{align}
\STATE \textbf{If} $M_{gap} \leq \varepsilon_{gap}$ hold, go to Step 5; \textbf{Else if} $r=R^{max}$, then quit;
\STATE \textbf{Else}
\STATE Call (\ref{Taylor}) with linearized point $\left(F_{mn,t}^r, \pi_{n,t}^r\right )$ and build the following problem \textbf{P2}:
\begin{align}
\label{SsocDRUC}
\begin{array}{*{20}{l}}
\textbf{P2}: \min ~(\ref{obj})+\varrho_{r+1}\sum_{t\in\mathcal{T}}\sum_{(m,n)\in \mathcal{L}_g}s_{mn,t}^{r+1}\\
\text{s.t.~(\ref{Rocof}), (\ref{fnadirMILP})-(\ref{BigMwindON}), (\ref{QSSlimit})}\\
\text{(\ref{GASbalance}), (\ref{avegasflow}), (\ref{Linepackpipe})-(\ref{Finallinepack}),  (\ref{compressors}), (\ref{GASsources}), (\ref{GASpress}), (\ref{coupling}), (\ref{WeymouthSOC}), (\ref{Taylor})}\\
\text{(\ref{logicrelation})-(\ref{Ramping}), (\ref{Rwind}), (\ref{scheduledWIND}), (\ref{LINElimit}), (\ref{Powerbalance}), (\ref{DRJCCtheorem})}.
\end{array}
\end{align}
\STATE Solve \textbf{P2} to obtain the solutions $F_{mn,t}^{r+1}$ and $\pi_{n,t}^{r+1}$. Calculate the maximum SOC relaxation gap $M_{gap}$.
\STATE Update  $\varrho_{r+1}:=\min\{\kappa \varrho_r, \varrho^{max}\}$, $r:=r+1$, and go to Step 2.
\STATE \textbf{End if}
\STATE Output the results.
\end{algorithmic}
\label{algo:SCP}
\end{algorithm}

We then develop a sequential algorithm based on PCCP \cite{Lipp-OE-2016} summarized in Algorithm 1 to solve (\ref{Taylor}). The convergence proof of Algorithm 1 can be found in \cite{Lipp-OE-2016}. In Algorithm 1, both \textbf{P1} and \textbf{P2} are MISOCP problems. The solutions of \textbf{P1} provide the initial linearized points for constraint (\ref{WeymouthCONCAVE}). The performance of Algorithm 1 is affected by the quality of initial point. To obtain a high-quality initial point, we add an additional linear penalty term $\rho\sum_{t\in\mathcal{T}}\sum_{(m,n)\in \mathcal{L}_g}(\pi_{m,t}-\pi_{n,t})$ (where $\rho$ is a predefined positive parameter) into the objective function of \textbf{P1} to tighten the SOC constraint (\ref{WeymouthSOC}).
\vspace{-0.1cm}
\section{Extension to including unimodality information}
In the proposed DR-FCUC model in (\ref{DRUC}), only mean and variance information is used for constructing the ambiguity set, which may lead to a conservative solution. In practice, wind power uncertainty is likely to be unimodal \cite{Li-TCNS-2019}. Incorporating unimodality information into the moment-based ambiguity set could reduce the solution conservativeness \cite{Shi-TSG-2019}, \cite{Li-TCNS-2019}.

We consider the DR joint chance constraint under the ambiguity set with both moment and unimodality information:
\begin{align}
\label{DRCCu}
\inf_{\mathbb{P}\in \mathcal{P}}\mathbb{P}\left\{P_{w,t}^W+R_{w,t}^W \leq \tilde{P}_{w,t}^{W}, \forall w \in \mathcal{W}\right\}\geq 1-\epsilon,\forall t\in \mathcal{T}
\end{align}
where the ambiguity set with moment and unimodality is
\begin{align}
\label{ambiguitysetu}
\mathcal{P}=\left\{\mathbb{P}:
\begin{array}{*{20}{l}}
\mathbb{E}_\mathbb{P}[\tilde{P}_{w,t}^{W}]=\mu_{w,t}, \mathbb{E}_\mathbb{P}\left[(\tilde{P}_{w,t}^{W}-\mu_{w,t})^2\right]=\sigma_{w,t}^2, \\
\tilde{P}_{w,t}^{W} \text{ is unimodal, } \forall w \in \mathcal{W}, t\in \mathcal{T}, \\
\end{array}
\right\}.
\end{align}

We reformulate the DR joint chance constraint (\ref{DRCCu}) as SOC constraints with using theorem 2.

\emph{Theorem 2: If $ \epsilon \leq 1/6$,  for all $t\in\mathcal{T}$, the DR joint chance constraint (\ref{DRCCu}) under the ambiguity set (\ref{ambiguitysetu}) is inner approximated by the following SOC constraints:}
\begin{subequations}
\label{DRJCCtheoremU}
\begin{align}
\label{UCHBF}
&r_{w,t}\sigma_{w,t}\leq -P_{w,t}^W-R_{w,t}^W+\mu_{w,t},\forall w \in \mathcal{W}\\
\label{soc1u}
&\left\|
\begin{array}{*{20}{l}}
2s_{w,t}\\
4/9
\end{array}\right\|_2 \leq 2\epsilon_{w,t}+4/9, \forall w\in\mathcal{W}\\
\label{soc2u}
&\left\|
\begin{array}{*{20}{l}}
4/3\\
r_{w,t}-s_{w,t}
\end{array}\right\|_2 \leq r_{w,t}+s_{w,t}, \forall w\in\mathcal{W}\\
\label{BFsumTu}
&\sum_{w\in\mathcal{W}}\epsilon_{w,t} \leq \epsilon\\
\label{Theoremvariableu}
&\epsilon_{w,t}, r_{w,t}, s_{w,t} \geq 0, \forall w \in \mathcal{W}.
\end{align}
\end{subequations}
\emph{Proof:} See the Appendix B.

The above derived results for the DR joint chance constraint under the ambiguity set with both moment and unimodaltity can be directly incorporated into the developed solution framework (i.e., Algorithm 1) by replacing (\ref{DRJCCtheorem}) with (\ref{DRJCCtheoremU}).

\vspace{-0.3cm}
\section{Case study}
We conduct numerical case studies on two test systems to demonstrate the effectiveness of the proposed approach. All programs are coded in MATLAB 2015a and solved by GUROBI via YALMIP toolbox \cite{Lofberg-C-2004} on a laptop with an inter(R) Core(TM) 2.6 GHz CPU and 16 GB memory. The MIPgap is set as 0.01. Parameters in Algorithm 1 are set as follows: $\varepsilon_{gap}=0.001$, $\varrho_0=0.02$, $\varrho^{\text{max}}=1000$, $\kappa=1.5$, and $R^{\text{max}}=50$. The system parameters are set as follows: load damping $D=1\%$ /Hz, frequency dead band $f_\text{DB}=15$ mHz, nominal frequency $f_0=50$ Hz, and delivery time $T_d=10$ s.
\vspace{-0.4cm}
\subsection{Sample generation and selection}
We use the Monte Carlo simulation to generate 20000 samples of wind power outputs from Gaussian distribution, whose mean is set to the wind power forecasted value and variance is set to 5\% of the mean value. We then divide the 20000 samples into two parts. The first part with 10000 samples forms the in-sample data. The second part is used in the out-of-sample test analysis. Assume that we have a limited knowledge on uncertain wind power $\tilde{P}_{w,t}^W$, we use a small number of samples from the in-sample data with $N=20$ to construct the ambiguity set with the empirical mean and variance. This setup is similar to \cite{Zhang-TPS-2017}.
\vspace{-0.4cm}
\subsection{The IEGS with a 5-bus system and a 7-node gas system}
A small IEGS with a 5-bus power system and a 7-node gas system is used here to test the performance of the proposed method. The test system consists of 1 non-GFU, 2 GFUs, 2 wind farms, 2 electricity loads, 2 gas sources, 1 compressor, and 3 gas loads. The total capacity of generators is 540 MW and the installed wind capacity is 200 MW (27\% wind penetration). The topology of the test system is illustrated in Fig.1. Detailed data of the test system is available in \cite{Data-2021}. The contingency event in this small IEGS is assumed to be a sudden increase of 5\% total electricity load. The frequency requirements are set as: $RoCoF^\text{max}=0.125$ Hz/s, $f^\text{min}=49.2$ Hz, and  $\triangle f^\text{max}_\text{qss}=0.2$ Hz \cite{Ding-TPS-2021}. The allowable joint violation probability $\epsilon$ is set as 0.05.

\begin{figure}[!t]
\centering
 \includegraphics[width=2.4in]{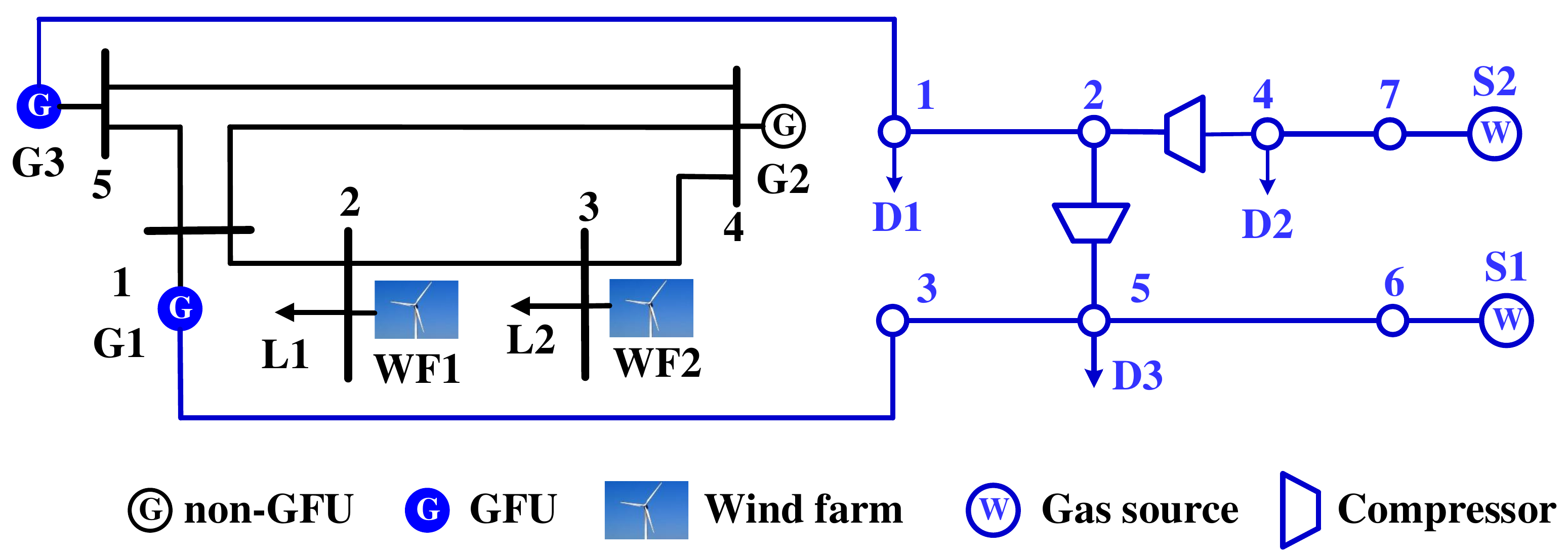}
\caption{A small IEGS with a 5-bus power system and a 7-node gas system.}
\label{E5_G7}
\end{figure}

\subsubsection{The Necessity of Including Frequency Constraints}
We compare the results of the following two scheduling models to verify the importance of incorporating frequency constraints (FCs).

i) Scheduling without FCs: the model is obtained from the DR-FCUC (\ref{DRUC}) by removing its FCs (\ref{Rocof}), (\ref{fnadirMILP})-(\ref{BigMwindON}), and (\ref{QSSlimit}) but adding capacity-based primary reserve constraints (i.e., $\sum_{i\in \mathcal{G} \cup \mathcal{N}}R_{i,t}^G+\sum_{w\in \mathcal{W}}R_{w,t}^W\geq \triangle P^{Loss}_i, \forall t \in \mathcal{T}$).

ii) Scheduling with FCs: the model is the propsoed DR-FCUC model in (\ref{DRUC}).

Table I reports the total cost, start-up and shut-down (SU/SD) cost, PFR cost of generators and WFSs, and VI provision cost of the scheduling models with/without FCs. Fig.2 and Fig.3 show the hourly post-contingency RoCoF and frequency nadir calculated by using the solutions of scheduling models with/without FCs.

As shown in Table I, with inclusion of FCs in the scheduling model, the SU/SD cost and PFR cost are reduced, while the generation cost and VI cost are increased. This could be explained as: the inclusion of FCs enforce that more generators and WFSs are committed to provide inertia. Therefore, the start-up/shut-down frequency of generators is reduced, leading to a less cost. Whereas generation cost and VI cost are increased. Overall, an evident increase in total cost is up to 1.5\% compared to the scheduling model without FCs. Although the scheduling model with FCs provides a more expensive solution, Fig.2 and Fig.3 show that both RoCoF and frequency nadir satisfy the frequency limits $RoCoF^{\text{max}}=0.125$ Hz/s and $f^{\text{min}}=49.2$ Hz in all hours. However, for scheduling without FCs, the RoCoF in hours 7-18, and 22-24 exceeds the maximum limit $RoCoF^{\text{max}}$ and the frequency nadir is below the minimum allowed frequency $f^{\text{min}}=49.2$ Hz in hours 16-17, 24. This observation indicates the scheduling without FCs would result in a potential frequency stability problem following a contingency.
\begin{table*}[!t]
\centering
\vspace{-0.3cm}
\label{fcs}
\caption{{Results of scheduling models with and without frequency constraints}}
\centering
\begin{tabular}{cccccc}
\hline
\hline
   Model &Total cost (\$) & SU/SD cost (\$) & Generation cost (\$) & PFR cost (\$)   &VI cost (\$) \\
\hline
 Scheduling without FC  &289726   &2510  &250886  &32730    &3600         \\
 Scheduling with FC      &293972   &1910   &252754  &32108  &7200         \\
\hline
\hline
\end{tabular}
\end{table*}
%
\begin{figure}[!t]
\centering
 \includegraphics[width=2.1in]{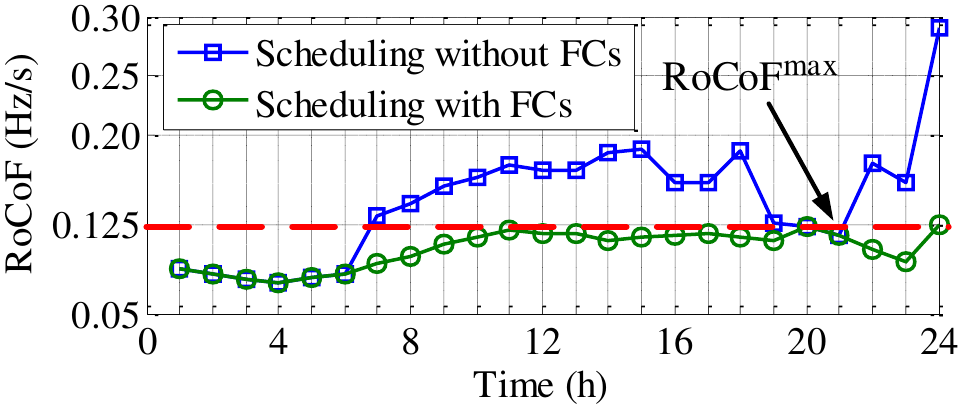}
\setlength{\abovecaptionskip}{-1.5pt}
\setlength{\belowcaptionskip}{-2.5pt}
\caption{RoCoF following the contingency at each hour produced by scheduling with/without FC.}
\label{RocoF}
\end{figure}
\begin{figure}[!t]
\centering
 \includegraphics[width=2.1in]{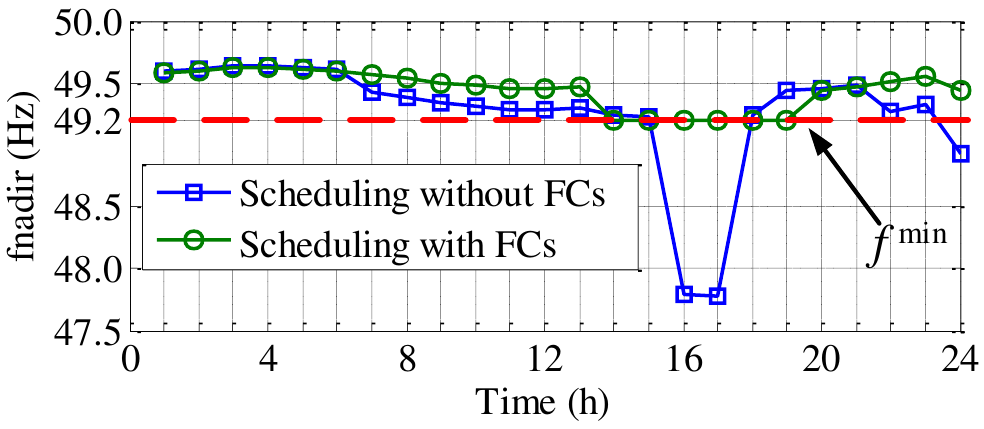}
\setlength{\abovecaptionskip}{-1.5pt}
\setlength{\belowcaptionskip}{-2.5pt}
\caption{Frequency nadir following the contingency at each hour produced by scheduling with/without FC.}
\label{Fnadir}
\end{figure}

\subsubsection{The Necessity of Including NGS Operational Constraints}
We compare the results of scheduling models with and without NGS operational constraints.

i) Scheduling with NGS operational constraints: the model is the proposed DR-FCUC model in (\ref{DRUC}).

ii) Scheduling without NGS operational constraints: the model is obtained from (\ref{DRUC}) by removing NGS operational constraints.

Note that we set a relatively low gas price in the test system to emphasise the power generation from GFUs.
Fig.4 and Fig.5 show the power outputs of generators without and with NGS operational constraints. For the scheduling without NGS operational constraints, the power generation from GFU G3 is largest in all hours, GFU G2 contributes to power generation at hours 8-21, while non-GFU G2 does not generate any power in most hours (i.e., 2-24). This is because, in the test system setting, the generation cost of GFUs (G1, G3) are cheaper than non-GFU (G2) and G3 is cheapest. With including NGS operational constraints, the power generation from GFU G3 decreases sharply while that of non-GFU G2 increases and is largest in most hours (i.e., 1-19). GFU G1  generates power in almost all hours and the generation curve is smoothed. Moreover, we find that including NGS operational constraints leads to a significant decrease on the total power generation from GFUs and thus results in higher total cost, as shown in Table II. These observations show that including NGS operational constraints have a marked impact on the operation of GFUs and the whole power system scheduling.

Furthermore, we fix variables associated with GFUs (i.e., $P_{g,t}$, $R_{g,t}^G$, $\forall g \in \mathcal{G}$) to the solutions obtained from the scheduling model without NGS operational constraints. We then build a gas network optimization problem with a constant objective function and NGS operational constraints (\ref{GASbalance})-(\ref{GASpress}). This problem is a nonlinear programming. We solve it by using IPOPT solver. However, the result shows that the gas network optimization problem is infeasible. This indicates that the solutions from the scheduling model without including NGS operational constraints are physically infeasible for the NGS and highlights the importance of including NGS operational constraints in the power system scheduling.
\begin{figure}[!t]
\centering
\vspace{-0.3cm}
 \includegraphics[width=2.1in]{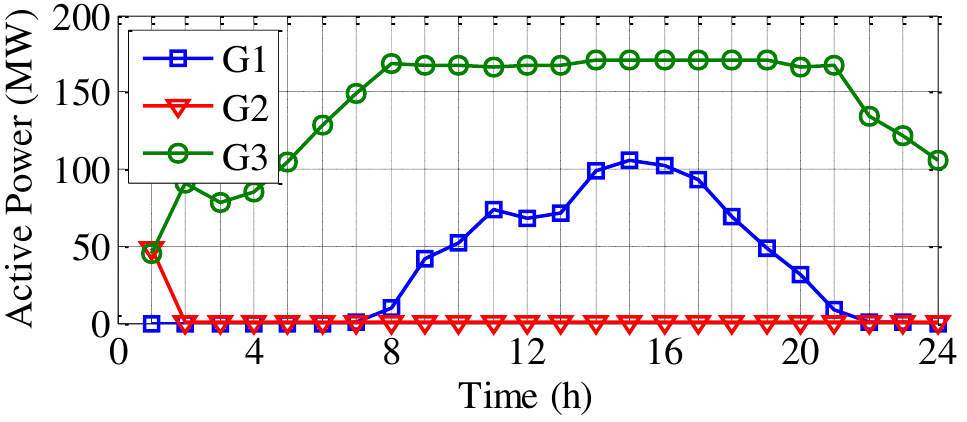}
\setlength{\abovecaptionskip}{-1.5pt}
\setlength{\belowcaptionskip}{-2.5pt}
\caption{Generation outputs without considering NGS operational constraints.}
\label{withoutGAS}
\end{figure}
\begin{figure}[!t]
\centering
\vspace{-0.3cm}
 \includegraphics[width=2.1in]{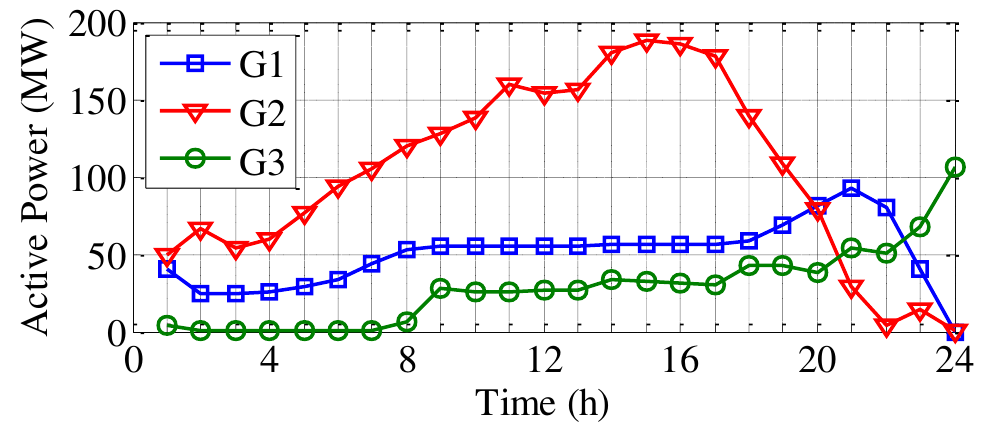}
\setlength{\abovecaptionskip}{-1.5pt}
\setlength{\belowcaptionskip}{-2.5pt}
\caption{Generation outputs with NGS operational constraints.}
\label{withGAS}
\end{figure}

\begin{table*}[t]
\centering
\vspace{-0.4cm}
\label{NGS}
\caption{{Results of scheduling models with and without NGS operational constraints}}
\centering
\begin{tabular}{cccccc}
\hline
\hline
   Model &Total cost (\$) & SU/SD cost (\$) & Generation cost (\$) & PFR cost (\$) & VI cost (\$) \\
\hline
  Scheduling without NGS operational constraints  &238925  &1910    &197738  &32108 &7200         \\
 Scheduling with NGS operational constraints      &293972   &1910   &252754  &32108  &7200        \\                                              \hline\hline
\end{tabular}
\end{table*}
\begin{table*}[!t]
\centering
\vspace{-0.4cm}
\label{VI}
\caption{{Results of scheduling models with and without VI provision from WFSs}}
\centering
\begin{tabular}{cccccc}
\hline
\hline
   Model &Total cost (\$) & SU/SD cost (\$) & Generation cost (\$) & PFR cost (\$) & VI cost (\$) \\
\hline
  Scheduling without VI   &296473   &1910   &251409  &43154 &0         \\
 Scheduling with VI       &293972   &1910   &252754  &32108  &7200        \\
\hline
\hline
\end{tabular}
\end{table*}

\subsubsection{Impacts of Virtual Inertia from WFSs}
We evaluate the influence of virtual inertia (VI) provision from wind farm systems (WFSs) on the system operation.
Table I reports the total cost, SU/SD cost, PFR cost of generators and WFSs, and VI provision cost of the scheduling models with and without VI from WFSs. We observe that the inclusion of VI provision from WFSs reduces the PFR cost but also introduces the VI cost. Overall, the total cost is reduced up to 0.8\% compared to the scheduling model without VI. This observation highlights the value of exploiting VI provision from WFSs in terms of operation cost saving.

\subsubsection{Comparisons with Other Chance Constrained FCUC}
We consider the following three FCUC models for comparisons.

i) SAA-FCUC: chance constrained FCUC model in \cite{Ding-TPS-2021}, which addresses the chance constraints by using sample average approximation. We term this model as SAA-FCUC.

ii) DR-FCUC-M: the proposed DR-FCUC model under the ambiguity set with \emph{moment} information only.

iii) DR-FCUC-U: the proposed DR-FCUC model under the ambiguity set with both moment and \emph{unimodality} information.

Note that SAA-FCUC model in \cite{Ding-TPS-2021} is dedicated to the power system. For fair comparison, we extend the SAA-FCUC model for an IEGS by adding the NGS operational constraints and then solve the model by Algorithm 1. We consider different sample sizes (10, 50, 100, 200, and 500) and compare the results: i) total cost, ii) empirical joint violation probability (EJVP) of solutions in out-of-sample test, and iii) computational efficiency. Note that the EJVP is defined as the percentage of samples for which any chance constraint is violated.

Table IV reports the total cost, EJVP, iterations, and computation time under different sample sizes for $\epsilon=0.05$. Compared SAA-FCUC with DR-FCUC-M/DR-FCUC-U, we see that the total costs of DR-FCUC-M and DR-FCUC-U are higher than that of SAA-FCUC, while their EJVPs are notably less than that of SAA-FCUC and all well-below the allowable joint violation probability $\epsilon=0.05$, indicating a higher solution reliability of DR-FCUC-M and DR-FCUC-U. On the contrary, although SAA-FCUC provides a solution with less cost, its EJVPs are well-above the allowable joint violation probability $\epsilon=0.05$ under all sample sizes. Although EJVPs of SAA-FCUC are decreased as the sample size grows, its computation time will grow dramatically as the increase of introduced binary variables, and the computation time with sample size being 500 is already up to 829.9 s. In comparison, the computation time consumed by DR-FCUC-M and DR-FCUC-U is far less than that of SAA-FCUC and within 30 s in this test system.  Furthermore, their computation time is not sensitive to the variation of sample size. In addition, we observe that these three models reuqire less than 5 iterations using Algorithm 1 under different sample sizes, showing a decent convergence performance of Algorithm 1.

Compared the results of DR-FCUC-M and DR-FCUC-U, we see that the total cost of DR-FCUC-U is less than that of DR-FCUC-M under different sample sizes, indicating the value of including unimodality information in the moment-based ambiguity set in terms of reducing solution conservativeness.

\begin{table}[!t]
 \centering
 \setlength{\abovecaptionskip}{0.0cm}   
 \caption{Comparisons with three FCUC models under different sample sizes for $\epsilon=0.05$}
 \begin{tabular}{cccccc}
 \hline\hline
      Model       &  $N$  & Total cost (\$)  & EJVP (\%)  & Iteration & Time (s)   \\
    \hline
                      &20   &262758  &93.00  &2 &18.8  \\

                      &50   &263444  &56.15  &2 &65.4\\
  SAA-FCUC            &100  &265631  &36.50  &3 &374.0  \\
                      &200  &268245  &20.50  &3  &476.3\\
                      &500  &270727  &11.60  &1  &829.9\\
    \hline
                      &20  &293972  &0  &2 &12.5    \\
                      &50  &295541  &0  &2 &8.7    \\
       DR-FCUC-M      &100 &295469 &0  &2 &11.9   \\
                      &200 &296367  &0  &2 &13.5   \\
                      &500 &296687  &0  &2 &14.5  \\
 \hline
                      &20  &278855  &1.38  &3  &17.3    \\
                      &50  &279974  &0.41   &3  &13.4   \\
       DR-FCUC-U      &100 &279932  &0.33   &3 & 14.9  \\
                      &200 &280577  &0.17  &3 &16.7   \\
                      &500 &280697  &0.15 &4 &20.9  \\
                      \hline\hline
\end{tabular}
\label{comparSIZE}
\end{table}

\vspace{-0.3cm}
\subsection{The IEGS with the IEEE 118-bus system and a 20-node gas system}
We use a larger IEGS with the IEEE 118-bus system and a 20-node gas system in \cite{Wang-TSE-2018} to demonstrate the scalability of the proposed methods.  The test system consists of 46 non-GFU, 8 GFUs, 6 wind farms, 99 electricity loads, 4 gas sources, 2 compressor, and 9 gas loads. The total capacity of generators is 9966.2 MW and the installed wind capacity is 3200 MW (24\% wind penetration). The topology and detailed data of the test system are available in \cite{Data-2021}. The contingency event considered is the loss of the largest generator with 805.2 MW. The frequency requirements in this test system are set as: $RoCoF^\text{max}=0.5$ Hz/s, $f^\text{min}=49.2$ Hz, and  $\triangle f^\text{max}_\text{qss}=0.2$ Hz. The allowable joint violation probability $\epsilon$ is set as 0.10. The sample size is 20.

We compare the SAA-FCUC, DR-FCUC-M, and DR-FCUC-U in terms of total cost, EJVP, iteration, and computation time. To verify the value of DR joint chance constraint, we also compare with the immediate models, which replace the DR joint chance constraints in DR-FCUC-M and DR-FCUC-U with DR \emph{individual} chance constraints. We term these two models as DR-FCUC-M-I and DR-FCUC-U-I. The allowable individual violation probability is set as 0.10. The results are summarized in Table V.

Similar to the results in Table IV, SAA-FCUC provides a solution with lower cost but extremely higher EJVP (well-above $\epsilon$) in comparison with that of DR-FCUC-M and DR-FCUC-U. Compared to DR-FCUC-M and DR-FCUC-U, we also observe that DR-FCUC-U offers a lower-cost but slightly higher-EJVP solution, indicating the value of including unimodality information into the ambiguity set in terms of reducing solution conservativeness. Additionally, the results in Table V show the computation time of DR-FCUC-M and DR-FCUC-U are significantly lower than that of SAA-FCUC and within 1 min.

Compared DR-FCUC-M/DR-FCUC-U with DR-FCUC-M-I/DR-FCUC-U-I, DR-FCUC-M-I/DR-FCUC-U-I provide a lower cost, but the calculated EJVPs are well-above the allowable violation probability 0.1. Therefore, DR joint chance constraint modeling provides a stronger guarantee on the solution reliability though it leads to a higher total cost in comparison with the DR individual chance constraint modeling.

\begin{table}[!t]
 \centering
 \setlength{\abovecaptionskip}{0.0cm}   
 \caption{Comparisons of the proposed DR-FCUC models and SAA-FCUC model for the larger IEGS under $\epsilon=0.10$}
 \begin{tabular}{ccccc}
 \hline\hline
      Model       & Total cost (\$)  & EJVP (\%)  & Iteration & Time (s)   \\
    \hline
  SAA-FCUC        &5531373  &99.88  &2 &441.5  \\
  DR-FCUC-M       &5825611  &0  &2  &42.7  \\
  DR-FCUC-U       &5685369  &0.41   &2 &51.5   \\
  DR-FCUC-M-I     &5589578  &39.50   &2 &27.2  \\
  DR-FCUC-U-I     &5533086  &99.47  &2 &45.1   \\
                      \hline\hline
\end{tabular}
\label{comparisonIEGS118}
\end{table}
\vspace{-0.2cm}
\section{Conclusion}
In this paper, we propose a DR-FCUC model considering IEGS operational constraints to co-optimize the unit commitment and virtual inertia of WFSs under wind power uncertainty. We formulate frequency constraints as MILP constraints and address nonconvex gas network constraints by using penalty convex-concave procedure. We separately derive SOC constraints for DR joint chance constraints under two types of ambiguity sets (i.e., moment information only and both moment and unimodality). Finally, the proposed DR-FCUC model is converted into an MISOCP. Numerical results demonstrate that the proposed DR-FCUC approach can provide a high-reliability and computationally efficient solution. Numerical results also show the necessity of including NGS operational constraints and the value of exploiting virtual inertia provision from WFSs in terms of cost saving, and verify that including unimodality information into the moment-based ambiguity set can lead to a less conservative solution.

Future work will derive FCs considering the detailed wind turbine dynamics to enhance the accuracy of the proposed DR-FCUC model.

\appendices
\section{Proof of Theorem 1}
First, the ambiguity set $\mathcal{D}$ defined in (\ref{ambiguityset}) satisfies Assumption (A1) in \cite{Xie-MP-2019}. Therefore, according to Theorem 3 in \cite{Xie-MP-2019}, the DR joint chance constraint (\ref{DRCC}) is equivalent to its Bonferroni approximation:
\begin{subequations}
\label{DRCCresult1}
\begin{align}
\label{DRICC}
\nonumber &\inf_{\mathbb{P}_{w,t}\in \mathcal{D}_{{w,t}}}\mathbb{P}_{w,t}\left\{P_{jw,t}^W+R_{w,t}^W \leq \tilde{P}_{w,t}^{W}\right\}\geq 1-\epsilon_{w,t}, \\ &\forall w \in \mathcal{W}, t\in \mathcal{T}\\
\label{BFsum}
&\sum_{w\in\mathcal{W}}\epsilon_{w,t} \leq \epsilon, \forall t\in \mathcal{T}\\
\label{BFpositive}
&\epsilon_{w,t}\geq 0, \forall w \in \mathcal{W}, t\in \mathcal{T}
\end{align}
\end{subequations}
where $\mathbb{P}_{w,t}$ and $\mathcal{D}_{w,t}$ are the probability distribution and the ambiguity set of $P_{w,t}^W$. The ambiguity set $\mathcal{D}_{w,t}$ is described as
  $\mathcal{D}_{w,t}=\left\{\mathbb{P}:\mathbb{E}_\mathbb{P}[\tilde{P}_{w,t}^{W}]=\mu_{w,t}, \mathbb{E}_\mathbb{P}\left[(\tilde{P}_{w,t}^{W}-\mu_{w,t})^2\right]=\sigma_{w,t}^2\right\} $.

Second, according to Theorem 3.1 in \cite{Calafiore-JOTA-2006}, the DR individual chance constraint (\ref{DRICC}) can be equivalently reformulated as
\begin{align}
\label{CHBFdrcc}
\sqrt{\frac{1-\epsilon_{w,t}}{\epsilon_{w,t}}}\sigma_{w,t}\leq -P_{w,t}^W-R_{w,t}^W+\mu_{w,t},\forall w \in \mathcal{W}, t \in \mathcal{T}.
\end{align}
Note that $\epsilon_{w,t}$ is an optimization variable, thus constraint (\ref{CHBFdrcc}) is nonconvex.

By introducing auxiliary variable $r_{w,t}$, constraint (\ref{CHBFdrcc}) is transformed into
\begin{subequations}
\begin{align}
\label{CHBF1}
&r_{w,t}\sigma_{w,t}\leq -P_{w,t}^W-R_{w,t}^W+\mu_{j,t},\forall w \in \mathcal{W}, t \in \mathcal{T}\\
\label{CHBF2}
&\sqrt{\frac{1-\epsilon_{w,t}}{\epsilon_{w,t}}} \leq r_{w,t}, \forall w \in \mathcal{W}, t \in \mathcal{T}.
\end{align}
\end{subequations}

Third, we observe $\sqrt{\frac{1-\epsilon_{w,t}}{\epsilon_{w,t}}}=\sqrt{\frac{1-\epsilon_{w,t}^2}{\epsilon_{w,t}(1+\epsilon_{w,t})}} \leq \sqrt{\frac{1}{\epsilon_{w,t}(1+\epsilon_{w,t})}}$ for $0 \leq \epsilon_{w,t}\leq \epsilon$. Then we propose the following conservative approximation for constraint (\ref{CHBF2}):
\begin{align}
\label{CHBF2approxi}
\sqrt{\frac{1}{\epsilon_{w,t}(1+\epsilon_{w,t})}} \leq r_{w,t}, \forall w \in \mathcal{W}, t \in \mathcal{T}.
\end{align}

By squaring both sides of constraint (\ref{CHBF2approxi}), we have
\begin{align}
\label{CHBF2approxi2}
\frac{1}{r_{w,t}^2} \leq \epsilon_{w,t}(1+\epsilon_{w,t}), \forall w \in \mathcal{W}, t \in \mathcal{T}.
\end{align}

By introducing auxiliary variable $s_{w,t}$, constraint (\ref{CHBF2approxi2}) is equvalently converted into
\begin{subequations}
\begin{align}
\label{TPsoc1}
&s_{w,t}^2 \leq \epsilon_{w,t}(1+\epsilon_{w,t}), \forall w \in \mathcal{W}, t \in \mathcal{T} \\
\label{TPsoc2}
&1 \leq s_{w,t}r_{w,t}, \forall w \in \mathcal{W}, t \in \mathcal{T}
\end{align}
\end{subequations}
whose standard SOC forms are (\ref{soc1}) and (\ref{soc2}), respectively.

Combining (\ref{BFsum}), (\ref{BFpositive}), (\ref{CHBF1}), (\ref{TPsoc1}), and (\ref{TPsoc2}),  we conclude the proof.

\section{Proof of Theorem 2}
The proof of Theorem 2 is similar to that of Theorem 1.

The ambiguity set $\mathcal{P}$ (\ref{ambiguitysetu}) also satisfies Assumption (A1) in \cite{Xie-MP-2019}. The DR joint chance constraint (\ref{DRCCu}) is equivalent to its Bonferroni approximation by Theorem 3 in \cite{Xie-MP-2019}:
\begin{subequations}
\label{DRCCresult2}
\begin{align}
\label{DRICCu}
\nonumber &\inf_{\mathbb{P}_{w,t}\in \mathcal{P}_{{w,t}}}\mathbb{P}_{w,t}\left\{P_{w,t}^W+R_{w,t}^W \leq \tilde{P}_{w,t}^{W}\right\}\geq 1-\epsilon_{w,t},\\&\forall w \in \mathcal{W}, t\in \mathcal{T}\\
\label{BFsumu}
&\sum_{w\in\mathcal{W}}\epsilon_{w,t} \leq \epsilon, \forall t\in \mathcal{T}\\
\label{BFpositiveu}
&\epsilon_{w,t}\geq 0, \forall w \in \mathcal{W}, t\in \mathcal{T}
\end{align}
\end{subequations}
where $\mathcal{P}_{w,t}=\left\{\mathbb{P}:\mathbb{E}_\mathbb{P}[\tilde{P}_{w,t}^{W}]=\mu_{w,t}, \mathbb{E}_\mathbb{P}\left[(\tilde{P}_{w,t}^{W}-\mu_{w,t})^2\right]\right.
\\ \left.=\sigma_{w,t}^2, \tilde{P}_{w,t}^{W} \text{~is unimodal}\right\} $.

By using the one-sided Vysochanskij-Petunin inequality \cite{Vysochanskij-1985}, for $0 \leq \epsilon \leq 1/6$, (\ref{DRICCu}) is reformulated as
\begin{align}
\label{CHBFdrccu}
\sqrt{\frac{4/9-\epsilon_{w,t}}{\epsilon_{w,t}}}\sigma_{w,t}\leq -P_{w,t}^W-R_{w,t}^W+\mu_{w,t},\forall w \in \mathcal{W}, t \in \mathcal{T}.
\end{align}

Constraint (\ref{CHBFdrccu}) is further transformed into
\begin{subequations}
\begin{align}
\label{CHBF1u}
&r_{w,t}\sigma_{w,t}\leq -P_{w,t}^W-R_{w,t}^W+\mu_{w,t},\forall w \in \mathcal{W}, t \in \mathcal{T}\\
\label{CHBF2u}
&\sqrt{\frac{4/9-\epsilon_{w,t}}{\epsilon_{w,t}}} \leq r_{w,t}, \forall w \in \mathcal{W}, t \in \mathcal{T}
\end{align}
\end{subequations}
where $r_{w,t}$ is an auxiliary variable.

Since $\sqrt{\frac{4/9-\epsilon_{w,t}}{\epsilon_{w,t}}}=\sqrt{\frac{16/81-\epsilon_{w,t}^2}{\epsilon_{w,t}(4/9+\epsilon_{w,t})}} \leq \sqrt{\frac{16/81}{\epsilon_{w,t}(4/9+\epsilon_{w,t})}}$ for $0 \leq \epsilon_{w,t}\leq \epsilon$, we propose the following conservative approximation for constraint (\ref{CHBF2u}):
\begin{align}
\label{CHBF2approxiu}
\sqrt{\frac{16/81}{\epsilon_{w,t}(4/9+\epsilon_{w,t})}} \leq r_{w,t}, \forall w \in \mathcal{W}, t \in \mathcal{T}.
\end{align}

By squaring both sides of constraint (\ref{CHBF2approxiu}), we have
\begin{align}
\label{CHBF2approxi2u}
\frac{16/81}{r_{w,t}^2} \leq \epsilon_{w,t}(4/9+\epsilon_{w,t}), \forall w \in \mathcal{W}, t \in \mathcal{T}.
\end{align}

By introducing auxiliary variable $s_{w,t}$, constraint (\ref{CHBF2approxi2u}) is equivalently converted into
\begin{subequations}
\begin{align}
\label{TPsoc1u}
&s_{w,t}^2 \leq \epsilon_{w,t}(4/9+\epsilon_{w,t}), \forall w \in \mathcal{W}, t \in \mathcal{T} \\
\label{TPsoc2u}
&4/9 \leq s_{w,t}r_{w,t}, \forall w \in \mathcal{W}, t \in \mathcal{T}
\end{align}
\end{subequations}
whose standard SOC forms are (\ref{soc1u}) and (\ref{soc2u}), respectively.

Constraints (\ref{BFsum}), (\ref{BFpositiveu}), (\ref{CHBF1u}), (\ref{TPsoc1u}), and (\ref{TPsoc2u}) are combined to conclude the proof.


\ifCLASSOPTIONcaptionsoff
  \newpage
\fi

\ifCLASSOPTIONcaptionsoff
  \newpage
\fi


\begin{thebibliography}{1}
\bibitem{Doherty-TPS-2010}
R. Doherty \emph{emph}, ``An assessment of the impact of wind generation on system frequency control," \emph{IEEE Trans. Power Syst.,} vol. 25, no. 1, pp. 452-460, Feb. 2010.
\bibitem{Xie-PIEEE-2011}
L. Xie \emph{et al.}, ``Wind integration in power systems: operational challenges and possible solutions," \emph{Proc. IEEE,} vol. 99, no. 1, pp. 214-232, Jan. 2011.
\bibitem{Matevosyan-Springer-2017}
J. Matevosyan and P. Du, ``Wind integration in ERCOT," in \emph{Integration of Large-Scale Renewable Energy into Bulk Power Systems}. Springer, 2017, pp. 1-25.
\bibitem{O'Sullivan-TPS-2014}
J. O'Sullivan \emph{et al.}, ``Studying the maximum instantaneous nonsynchronous generation in an island system-frequency stability challenges in Ireland,” \emph{IEEE Trans. Power Syst.,} vol. 29, no. 6, pp. 2943-2951, Nov. 2014.
\bibitem{Restrepo-TPS-2005}
J. F. Restrepo and F. D. Galiana, ``Unit commitment with primary frequency regulation constraints," \emph{IEEE Trans. Power Syst.,} vol. 20, no. 4, pp. 1836-1842, Nov. 2005.
\bibitem{Chang-TPS-2013}
G. W. Chang, C. Chuang, T. Lu, and C. Wu, ``Frequency-regulating reserve constrained unit commitment for an isolated power system," \emph{ IEEE Trans. Power Syst.,} vol. 28, no. 2, pp. 578–586, May 2013.
\bibitem{Ahmadi-TPS-2014}
H. Ahmadi and H. Ghasemi, ``Security-constrained unit commitment with linearized system frequency limit constraints," \emph{IEEE Trans. Power Syst.,} vol. 29, no. 4, pp. 1536-1545, Jul. 2014.
\bibitem{Zhang-TPS-2020}
Z. Zhang, E. Du, F. Teng, N. Zhang and C. Kang, ``Modeling frequency dynamics in unit commitment with a high share of renewable energy," \emph{IEEE Trans. Power Syst.,} vol. 35, no. 6, pp. 4383-4395, Nov. 2020.
\bibitem{Trovato-TPS-2019}
V. Trovato, A. Bialecki and A. Dallagi, ``Unit commitment with inertia-dependent and multispeed allocation of frequency response services," \emph{IEEE Trans. Power Syst.,} vol. 34, no. 2, pp. 1537-1548, Mar. 2019.
\bibitem{Lee-TPS-2013}
Y. Y. Lee and R. Baldick, ``A frequency-constrained stochastic economic dispatch model," \emph{IEEE Trans. Power Syst.,} vol. 28, no. 3, pp. 2301-2312, Aug. 2013.
\bibitem{Teng-TPS-2016}
F. Teng, V. Trovato and G. Strbac, ``Stochastic scheduling with inertia-dependent fast frequency response requirements," \emph{IEEE Trans. Power Syst.,} vol. 31, no. 2, pp. 1557-1566, Mar. 2016.
\bibitem{Paturet-TPS-2020}
M. Paturet, U. Markovic, S. Delikaraoglou, E. Vrettos, P. Aristidou and G. Hug, ``Stochastic unit commitment in low-inertia grids," \emph{IEEE Trans. Power Syst.,} vol. 35, no. 5, pp. 3448-3458, Sept. 2020.
\bibitem{Wen-TPS-2016}
Y. Wen, W. Li, G. Huang, and X. Liu, ``Frequency dynamics constrained unit commitment with battery energy storage," \emph{IEEE Trans. Power Syst.,} vol. 31, no. 6, pp. 5115-5125, Nov. 2016.
\bibitem{Prakash-TPS-2018}
V. Prakash, K. C. Sharma, R. Bhakar, H. P. Tiwari and F. Li, ``Frequency response constrained modified interval scheduling under wind uncertainty," \emph{IEEE Trans. Sustain. Energy,} vol. 9, no. 1, pp. 302-310, Jan. 2018.
\bibitem{Ding-TPS-2021}
T. Ding, Z. Zeng, M. Qu, J. P. S. Catal\~{a}o, and M. Shahidehpour, ``Two-stage chance-constrained stochastic unit commitment for optimal provision of virtual inertia in wind-storage systems," \emph{IEEE Trans. Power Syst.,} vol. 36, no. 4, pp. 3520-3530, Jul. 2021.
\bibitem{Wen-TSG-2018}
Y. Wen, X. Qu, W. Li, X. Liu, and X. Ye, ``Synergistic operation of electricity and natural gas networks via ADMM," \emph{IEEE Trans. Smart Grid,} vol. 9, no. 5, pp. 4555-4565, Sep. 2018.
\bibitem{Chen-TPS-2020OE}
S. Chen, A. J. Conejo, R. Sioshansi and Z. Wei, ``Operational equilibria of electric and natural gas systems with limited information interchange,"  \emph{IEEE Trans. Power Syst.,} vol. 35, no. 1, pp. 662-671, Jan. 2020.
\bibitem{Fang-TSE-2018}
J. Fang, Q. Zeng, X. Ai, Z. Chen and J. Wen, ``Dynamic optimal energy flow in the integrated natural gas and electrical power systems,"  \emph{IEEE Trans. Sustain. Energy,} vol. 9, no. 1, pp. 188-198, Jan. 2018.
\bibitem{Wang-TSE-2018}
C. Wang, W. Wei, J. Wang, L. Bai, Y. Liang, and T. Bi, ``Convex  optimization based distributed optimal gas-power flow calculation”,
\emph{IEEE Trans. Sustain. Energy,} vol. 9, no. 3, pp. 1145-1156, Jul. 2018.
\bibitem{Zhao-TPS-2018}
B. Zhao, A. J. Conejo, and R. Sioshansi, ``Coordinated expansion planning of natural gas and electric power systems," \emph{IEEE Trans. Power Syst.,} vol. 33, no. 3, pp. 3064-3075, May 2018.
\bibitem{Chen-TPS-2020}
S. Chen, A. J. Conejo, R. Sioshansi and Z. Wei, ``Equilibria in electricity and natural gas markets with strategic offers and bids" \emph{IEEE Trans. Power Syst.,} vol. 35, no. 3, pp. 1956-1966, May 2020.
\bibitem{Zheng-TSG-2021}
W. Zheng, W. Huang, D. J. Hill and Y. Hou, ``An adaptive distributionally robust model for three-Phase distribution network reconfiguration," \emph{IEEE Trans. Smart Grid,} vol. 12, no. 2, pp. 1224-1237, Mar. 2021.
\bibitem{Zhang-TPS-2017}
Y. Zhang, S. Shen, and J. L. Mathieu, ``Distributionally robust chance constrained optimal power flow with uncertain renewables and uncertain reserve provided by loads,"  \emph{IEEE Trans. Power Syst.}, vol. 32, no. 2, pp, 1378-1388, Mar. 2017.
\bibitem{Zare-TPS-2018}
A. Zare, C. Y. Chung, J. Zhan and S. O. Faried, ``A distributionally robust chance-constrained MILP model for multistage distribution system planning with uncertain renewables and loads," \emph{IEEE Trans. Power Syst.,} vol. 33, no. 5, pp. 5248-5262, Sept. 2018.
\bibitem{Shi-TSG-2019}
Z. Shi, H. Liang, S. Huang, and V. Dinavahi, ``Distributionally robust chance-constrained energy management for islanded microgrids
," \emph{IEEE Trans. Smart Grid}, vol. 10, no. 2, pp. 2234-2244, Mar. 2019.
\bibitem{Fang-APEN-2019}
X. Fang, B. M. Hodge, H. Jiang, and Y. Zhang, ``Decentralized wind uncertainty management: Alternating direction method of multipliers based distributionally-robust chance constrained optimal power flow," \emph{Appl. Energy}, vol. 239, no. 1, pp, 938-947, Apr. 2019.
\bibitem{Yang-TSG-2021}
L. Yang, Y. Xu, W. Gu and H. Sun, ``Distributionally robust chance-constrained optimal power-gas flow under bidirectional interactions considering uncertain wind power," \emph{IEEE Trans. Smart Grid,}, vol. 12, no. 2, pp. 1722-1735, March 2021.
\bibitem{Chu-TSG-2021}
Z. Chu, N. Zhang and F. Teng, ``Frequency-constrained resilient scheduling of microgrid: A distributionally robust approach," \emph{IEEE Trans. Smart Grid}, 2021. doi: 10.1109/TSG.2021.3095363.
\bibitem{Lipp-OE-2016}
T. Lipp, and S. Boyd, ``Variations and extension of the convex-concave procedure," \emph{Optim. Eng.}, vol. 17, pp. 263-287, 2016.
\bibitem{Li-TCNS-2019}
B. Li, R. Jiang, and J. L. Mathieu, ``Distributionally robust chance-constrained optimal power flow assuming unimodal distributions with misspecified modes," \emph{IEEE Trans. Control. Network Syst.}, vol. 6, no. 3, pp. 1223-1234, Sep. 2018.
\bibitem{Lofberg-C-2004}
J. Lofberg, ``YALMIP: A toolbox for modeling and optimization in MATLAB," \emph{IEEE International Conference on Robotics and Automation,} 2004, pp. 284-289.
\bibitem{Data-2021}
\emph{Topologies and data of two IEGSs used for case studies}, [Online] Available:
https://cloud.tsinghua.edu.cn/f/d44d7f82222f4c2a8592/
\bibitem{Xie-MP-2019}
W. Xie, S. Ahmed, and R. Jiang, ``Optimized Bonferroni approximations of distributionally robust joint chance constraints," \emph{Math. Program.}, Nov. 2019. https://doi.org/10.1007/s10107-019-01442-8.
\bibitem{Calafiore-JOTA-2006}
G. C. Calafiore and L. El Ghaoui, ``On distributionally robust chance constrained linear programs," \emph{J. Optim. Theory Appl.}, vol. 130, no. 1, pp. 1-22, 2006.
\bibitem{Vysochanskij-1985}
D. Vysochanskij and Y. Petunin, ``Improvement of the unilateral 3 $\sigma$ -rule for unimodal distributions," Dokl. Akad. Nauk. Ukr. SSR, Ser. A, vol. 1, pp. 6-8, 1985.
\end{thebibliography}
\end{document}